\documentclass[12pt]{article}
\usepackage{epsfig,amsmath,amssymb,stmaryrd,enumerate}
\usepackage{graphicx}

\advance \textheight by \topskip
\advance \textheight by 1pt
\makeatother

\begin{document}

\title{Integrable geodesic flows on 2-torus: formal solutions and variational principle}

\author{{Misha Bialy\thanks{School of Mathematical Sciences, Raymond and Beverly Sackler Faculty of Exact Sciences,
Tel Aviv University, Israel; e-mail: bialy@post.tau.ac.il} \ and
Andrey E. Mironov\thanks{Sobolev Institute of Mathematics and
Novosibirsk State University, Russia; e-mail: mironov@math.nsc.ru}}}

\date{}
\maketitle

\begin{abstract}In this paper we study quasi-linear system of
partial differential equations which describes the existence of  the
polynomial in momenta first integral of the integrable geodesic flow
on 2-torus. We proved in \cite{BM} that this is a semi-Hamiltonian system
and we show here that the metric associated with the system is a
metric of Egorov type. We use this fact in order to prove that in
the case of integrals of degree three and four the system is in fact
equivalent to a single remarkable equation of order 3 and 4
respectively. Remarkably the equation for the case of degree four
has variational meaning: it is Euler-Lagrange equation of a
variational principle. Next we prove that this equation for $n=4$
has formal double periodic solutions as a series in a small
parameter.

\end{abstract}

\section{Introduction}
In this paper we study integrable geodesic flows on two-dimensional
torus $\mathbb{T}^2=\mathbb{R}^2/\Gamma$, where
$\Gamma\subset\mathbb{Z}^2$ is a lattice. Let
$ds^2=\sum_{i,j=1}^2g_{ij}(q)dq^idq^j$ be a Riemannian metric on
2-torus. The geodesic flow of the metric is called integrable if the
Hamiltonian system
$$
 {\dot q}^i=\frac{\partial H}{\partial p_i},\quad \dot{p}^i=-\frac{\partial H}{\partial q^i},\quad H=\frac{1}{2}\sum_{i,j=1}^2g^{ij}(q)p_ip_j,
$$
has a first integral $F(q,p):T^*\mathbb{T}^2\rightarrow {\mathbb
R}$, i.e.
$$
 \dot{F}=\{F,H\}=\left(\frac{\partial H}{\partial q^1}\frac{\partial F}{\partial p_1}-\frac{\partial H}{\partial p_1}\frac{\partial F}{\partial q^1}\right)+
 \left(\frac{\partial H}{\partial q^2}\frac{\partial F}{\partial p_2}-\frac{\partial H}{\partial p_2}\frac{\partial F}{\partial q^2}\right)=0,
$$
such that almost everywhere $F$ is independent with $H$.
There are two kind of Riemannian metrics having polynomial
integrals for the geodesic flow. If the metric has the form $$ds^2=\Lambda(\alpha x+\beta y)(dx^2+dy^2)\
or \ ds^2=(\Lambda_1(\alpha_1 x+\beta_1 y)+\Lambda_2(\alpha_2 x+\beta_2 y))(dx^2+dy^2),$$ then there exist polynomial integrals of
 degree one or degree two.
The existence of Riemannian metrics with  non-reducible polynomial integrals of higher degree is a difficult  open problem.

Let us recall some results on integrable geodesic flows on
$\mathbb{T}^2$. If the geodesic flow is integrable, then on the
torus there are global semi-geodesic coordinates $(t,x)$ (see
\cite{BM}), i.e. $$ds^2=g^2(t,x)dt^2+dx^2,\
H=\frac{1}{2}\left(\frac{p_1^2}{g^2}+p_2^2\right).$$ The polynomial
first integral has the form
$$
 F=\frac{a_0}{g^n}p_1^n+\frac{a_1}{g^{n-1}}p_1^{n-1}+\dots+\frac{a_{n-2}}{g}p_1^2p_2^{n-2}+p_1p_2^{n-1}+p_2^n,\quad a_s=a_s(t,x).
$$
The condition $\{F,H\}=0$ is equivalent to the quasi-linear system of partial differential equations
\begin{equation}\label{e1}
 U_t+A(U)U_x=0,
\end{equation}
where $U^{\top}=(a_0,\dots,a_{n-1})$, $a_{n-1}=g$,
$$
 A=  \left(
  \begin{array}{cccccc}
   0 & 0 & \dots &
  0 & 0 & a_1  \\
  a_{n-1} &
 0 & \dots & 0 & 0 & 2a_2-na_0\\0 &
 a_{n-1} & \dots & 0 & 0 & 3a_3-(n-1)a_1\\
 \dots & \dots & \dots & \dots & \dots & \dots \\
 0 &
 0 & \dots & a_{n-1} & 0 & (n-1)a_{n-1}-3a_{n-3}\\
 0 &
 0 & \dots & 0 & a_{n-1} & na_n-2a_{n-2}\\
  \end{array}\right).
$$
System (\ref{e1}) has remarkable properties. It can be written
in the form of conservation laws, i.e. there is a regular change of variables
$$
 (a_0,\dots,a_{n-1})\mapsto (G_1(a),\dots,G_n(a))
$$
such that for some $F_1(a),\dots,F_n(a)$ the conservation laws
hold true
$$
 (G_i(a))_t+(F_i(a))_x=0,\quad i=1,\dots,n.
$$

Moreover, in the hyperbolic domain, where eigenvalues $\lambda_1,\dots,\lambda_n$ of $A$ are real and pairwise distinct there exist change of variable
$$
 (a_0,\dots,a_{n-1})\mapsto (r_1(a),\dots,r_n(a))
$$
such that the system can be written in Riemann invariants.
$$
 (r_i)_t+\lambda_i(r)(r_i)_x=0,\quad i=1,\dots,n.
$$
Such systems are called semi-Hamiltonian and generalized hodograph method applies \cite{Ts}. Semi-Hamiltonian systems and in particular
systems of hydrodynamic type are very important in mathematical physics (see for example \cite{DN}, \cite{MF}, \cite{S}).

 It is not clear to us at present if and how the generalized hodograph method can be used to prove non-existence of
smooth solutions for a semi-Hamiltonian system. Using the original
idea by P. Lax based on analysis along characteristics we proved for
the cases $n=3,4$ that in the elliptic domain (where matrix $A$ has
two complex-conjugated eigenvalues) the behavior of solutions can be
analyzed: the integrals of degree three and four are reduced to
integrals of degree one or two \cite{BM1}. Thus nontrivial integrals
of degree $3,4$ may exist only in the hyperbolic region of the
quasi-linear system (\ref{e1}).

In this paper we study a quasi-linear system (\ref{e3}) which
corresponds to the choice of conformal coordinates $(x,y)$ for Riemannian
metric $ds^2=\Lambda(dx^2+dy^2)$.
  We assume that the geodesic flow has a polynomial in momenta integral
$$
 F=a_0(x,y)p_1^n+a_1(x,y)p_1^{n-1}p_2+\dots+a_n(x,y)p_2^n.
$$
Kozlov and Denisova  \cite{KD} proved that if $\Lambda$ is trigonometric polynomial then the
geodesic flow has no irreducible polynomial integrals of degree higher than two (see also \cite{KT}).
By Kolokoltsov's \cite{K} theorem
$$
 a_{n-1}=c_1+a_{n-3}-a_{n-5}+\dots, \quad a_n=c_2+a_{n-2}-a_{n-4}+\dots,
$$
where $c_1,c_2$ are some constants. Then the condition
$\{H,F\}=0,$ where $H=\frac{p_1^2+p_2^2}{2\Lambda}$ is equivalent to the system of quasi-linear equations
\begin{equation}\label{e3}
 A(U)U_x+B(U)U_y=0,\quad U=(a_0,a_1,\dots,a_{n-2},\Lambda)
\end{equation}
(see ({\ref{eq8}) below to specify $A(U),B(U)$ explicitly).
This system also can be written in the form of conservation laws and moreover, in the
hyperbolic region it admits $n$ Riemann invariants, so the system is
semi-Hamiltonian (see \cite {b}).

 Let us remind that for semi-Hamiltonian systems the following relations on eigenvalues hold
$$
 \partial_{r_j}\frac{\partial_{r_i}\lambda_k}{\lambda_i-\lambda_k}=\partial_{r_i}
 \frac{\partial_{r_j}\lambda_k}{\lambda_j-\lambda_k},\quad i\ne j\ne k\ne i.
$$
These relations mean that there exists a diagonal metric on the space of field variables
\begin{equation}\label{e2}
 ds^2=H_1^2(r)dr_1^2+\dots+H_n^2(r)dr_n^2
\end{equation}
with Christoffel symbols satisfying the identities
$$
 \Gamma_{ki}^k=\frac{\partial_{r_i}\lambda_k}{\lambda_i-\lambda_k},\quad i\ne k.
$$

Let us formulate now our main results. In Theorems 1--7 we assume that $c_1=0.$ This can be achieved by a rotation in the plane $x,y$.

\vspace{0.4cm}

\noindent {\bf Theorem 1} {\it The metric (\ref{e2}) associated with
the semi-Hamiltonian system (\ref{e3}) is a metric of Egorov type,
i.e. the rotation coefficients
$$
 \beta_{ij}=\frac{\partial_{r_i}H_j}{H_i},\quad i\ne j
$$
are symmetric $\beta_{ij}=\beta_{ji}$, or equivalently there is a function $A(r)$ such that
$$
 \partial_{r_i}A(r)=H_i^2(r).
$$}

\vspace{0.4cm}

\noindent In fact it follows from theorem of Pavlov and Tsarev
\cite{PT} that in order to prove Theorem 1 one needs to find two
conservation laws of a special form. In the next theorem we state
 the existence of these conservation laws for (\ref{e3}).

\vspace{0.4cm}

\noindent {\bf Theorem 2} {\it The system (\ref{e3}) has two
conservation laws of the form
\begin{equation}\label{eqq3}
 P(U)_x+Q(U)_y=0,\quad Q(U)_x+R(U)_y=0.
\end{equation}
}

\vspace{0.4cm}

\noindent Functions $P(U),Q(U),R(U)$ are found explicitly in Lemma 1
(see below). Remarkably, for $n=3,4$  Theorem 2 allows us to reduce
system (\ref{e3}) to a single equation.

\vspace{0.4cm}

\noindent{\bf Theorem  3} {\it Let $n=3$, and $\lambda(x,y)$ be a
solution periodic with respect to the lattice $\Gamma$  of the
equation
\begin{equation}\label{e35}
 \Delta\lambda=\frac{3c_2}{2}\Lambda-2a_{11}-2a_{22}.
\end{equation}
Then the function $\lambda$ satisfies the equation
$$
2\lambda_{xx}\lambda_{xxy}+\lambda_{yy}(\lambda_{xxy}-\lambda_{yyy})+\lambda_{xy}(\lambda_{xxx}+\lambda_{xyy})+4a_{11}\lambda_{xxy}+
$$
\begin{equation}\label{eq4}
 2a_{12}(\lambda_{xxx}+\lambda_{xyy})+
2a_{22}(\lambda_{xxy}-\lambda_{yyy})=0,
\end{equation}
where $a_{11},a_{12},a_{22}$ are some constants defined by the metric and the integral.  }

\vspace{0.4cm}

\noindent{\bf Theorem  4} {\it Let $n=4$, and $\lambda(x,y)$ be a
solution periodic with respect to the lattice $\Gamma$  of the
equation
$$
 \Delta\lambda=2c_2\Lambda-2a_{11}-2a_{22}.
$$
Then the function $\lambda$ satisfies the equation
$$
 \lambda_{xy}(\lambda_{yyyy}-\lambda_{xxxx})+
 3(\lambda_{yyy}\lambda_{xyy}-\lambda_{xxy}\lambda_{xxx})
+
$$
\begin{equation}\label{e4}
+2(\lambda_{yy}\lambda_{xyyy}-\lambda_{xx}\lambda_{xxxy})+4a_{22}\lambda_{xyyy}-4a_{11}\lambda_{xxxy}+
2 a_{12}(\lambda_{yyyy}-\lambda_{xxxx})=0,
\end{equation}
where $a_{11},a_{12},a_{22}$ are some constants defined by the metric and the integral.  }

\vspace{0.4cm}

Our next theorem states that the equations (\ref{eq4}) and
(\ref{e4}) are in fact equivalent to the system (\ref{e3}).

\vspace{0.4cm}

\noindent{\bf Theorem  5} {\it Let $\lambda$ be a solution of the
equation (\ref{eq4}) or (\ref{e4}) respectively periodic with
respect to the lattice $\Gamma$
 which satisfies the condition
$\Delta\lambda+2(a_{11}+a_{22})>0.$ Then the corresponding
 solution of the system (\ref{e3}) is periodic also. }

\vspace{0.4cm}

\noindent Proof of
this theorem is very simple for the case $n=3$ but requires certain
topological argument for the case $n=4$. This argument is given
below in Section 4.

Remarkably equation (\ref{e4}) admits the following variational
interpretation.

\vspace{0.4cm}

\noindent{\bf Theorem  6} {\it Equation (\ref{e4}) coincides with
the Euler-Lagrange equation of the functional

$$\mathcal{L}(\lambda)=\int \frac{1}{2}\left (
4\lambda_{xy}(a_{22}\lambda_{yy}-a_{11}\lambda_{xx})+2a_{12}(\lambda_{yy}^2-\lambda_{xx}^2)
+\lambda_{xy}(\lambda_{yy}^2-\lambda_{xx}^2) \right ) dx\ dy$$ }

\vspace{0.4cm}

Let us remark that the functional becomes especially simple
$$\mathcal{L}(\lambda)=\int \frac{1}{2}\lambda_{xy}(\lambda_{yy}-\lambda_{xx})(\frac{1}{\varepsilon}+
\Delta\lambda) dx\ dy
$$ for the choice of constants $a_{ij}$ where
$a_{11}=a_{22}=\frac{1}{4\varepsilon}, a_{12}=0.$ However, at the present moment we have
no significant results on the critical points of this functional.

 \vspace{0.4cm}

 \vspace{0.4cm}
Let us consider now in more details the equation (\ref{e4}) for
$a_{12}=0, a_{11}=a_{22}=\frac{1}{4\varepsilon}$
$$
 \lambda_{xxxy}-\lambda_{xyyy}=\varepsilon(\lambda_{xy}(\lambda_{yyyy}-\lambda_{xxxx})+
 3(\lambda_{yyy}\lambda_{xyy}-\lambda_{xxy}\lambda_{xxx})
+
$$
\begin{equation}\label{e5}
 +2(\lambda_{yy}\lambda_{xyyy}-\lambda_{xx}\lambda_{xxxy})).
\end{equation}
The conformal factor of the metric has the form
$$
 \Lambda=\frac{1}{c_2}(\frac{\Delta\lambda}{2}+\frac{1}{2\varepsilon}).
$$
Let us look for a solution of (\ref{e5}) as a formal power series in
$\varepsilon$:
\begin{equation}\label{es5}
 \lambda(x,y)=\lambda_0(x,y)+\lambda_1(x,y)\varepsilon+\lambda_2(x,y)\varepsilon^2+\dots,
\end{equation}
where $\varepsilon$ is a small parameter. Then from (\ref{e5}) we
have a recursion formula.
\begin{equation}\label{e6}
 (\lambda_k)_{xyyy}-(\lambda_k)_{xxxy}=\sum_{s=0}^{k-1}<\lambda_s,\lambda_{k-s-1}>
\end{equation}
where
$$
 <\lambda_p,\lambda_q>=\lambda_{p_{xy}}(\lambda_{q_{yyyy}}-\lambda_{q_{xxxx}})+3(\lambda_{p_{yyy}}\lambda_{q_{xyy}}-
 \lambda_{p_{xxx}}\lambda_{q_{xxy}})+
 2(\lambda_{p_{yy}}\lambda_{q_{xyyy}}-\lambda_{p_{xx}}\lambda_{q_{xxxy}}).
$$
In the rest of this section we shall assume that
$\Gamma\subset\mathbb{R}^2$ is the integer lattice $\mathbb{Z}^2$.
Given initial doubly periodic function $\lambda_0$ we wish to find
all $\lambda_k$ recursively by means of equation (\ref{e6}). It is
easy to see that (\ref{e6}) has a periodic solution $\lambda_k$ if
the Fourier series of the right hand side does not have monomials of
the form
$$
 e^{inx},\quad e^{iny},\quad e^{in(x+y)},\quad e^{in(x-y)}.
$$ Let us remark  also that the periodic
solution $\lambda_k$ of (\ref{e6}) is defined up to addition to
$\lambda_k$ some function of the form
$$
 \tilde{\lambda}_k=f_1(x)+f_2(y)+f_3(x-y)+f_4(x+y).
$$

Next theorem states that recursive process is well defined and thus
gives formal periodic solution of (\ref{e5}) if the initial function
$\lambda_0$ and the additions $\tilde{\lambda}_k$ on every stage $k$
are chosen symmetric with respect to coordinate axes and diagonals,
$x=0, y=0, x=y, x=-y.$ It is an important open question if the
convergence of the series can be achieved by a good choice of
initial function $\lambda_0$ and the functions $\tilde{\lambda}_k.$

\vspace{0.4cm}

\noindent{\bf Theorem  7} {\it Let
$$
 \lambda_0=\sum_{n\in{\mathbb N}}\alpha_{n}(\cos(nx)+\cos(ny))+
 \sum_{n\in{\mathbb N}}\beta_{n}(\cos(n(x-y))+\cos(n(x+y))),
$$
$$
 \tilde{\lambda}_k=\sum_{n\in{\mathbb N}}\alpha^k_{n}(\cos(nx)+\cos(ny))+
 \sum_{n\in{\mathbb N}}\beta^k_{n}(\cos(n(x-y))+\cos(n(x+y))),
$$
then the recursion formula gives a well defined formal periodic
solution (\ref{es5}) of (\ref{e6}).}

\vspace{0.4cm}

We prove Theorem 7 in the last section.

\vspace{0.4cm}

\noindent{\bf Remark}
{\it It is an interesting fact that substituting into the equation (\ref{e5})
 $$
 \lambda=f_1(x)+f_2(y)+f_3(x-y)+f_4(x+y)
$$ one gets an equation on the functions $f_1,...,f_4$ which was studied recently in
\cite{AA}. It was proved later in \cite{KTD} that there are no new
periodic solutions of this equation in this particular form. }


\section{Proof of Theorems 1 and 2}
Let us assume that $ds^2=\Lambda(x,y)(dx^2+dy^2)$ is a metric on
$\mathbb{T}^2$ and $F$ is a first integral polynomial in momenta for
the geodesic flow
$$
 F=a_0(x,y)p_1^n+a_1(x,y)p_1^{n-1}p_2+\dots+a_n(x,y)p_2^n.
$$
We have
$$
 2\Lambda^2\{H,F\}=2 (a_{0_x}p_1^n+a_{1_x}p_1^{n-1}p_2+a_{2_x}p_1^{n-2}p_2^2+\dots+a_{n_x}p_2^n)p_1\Lambda+
$$
$$
 +(na_0p_1^{n-1}+(n-1)a_1p_1^{n-2}p_2+(n-2)a_2p_1^{n-3}p_2^2+\dots+a_{n-1}p_2^{n-1})(p_1^2+p_2^2)\Lambda_x+
$$
$$
 2(a_{0_y}p_1^n+a_{1_y}p_1^{n-1}p_2+a_{2_y}p_1^{n-2}p_2^2+\dots+a_{n_y}p_2^n)p_2\Lambda+
$$
$$
 +(a_1p_1^{n-1}+2a_2p_1^{n-2}p_2+3a_3p_1^{n-3}p_2^2+\dots+na_{n}p_2^{n-1})(p_1^2+p_2^2)\Lambda_y.
$$
Let us denote by $l_{n+1},\dots,l_0$ coefficients of the homogeneous polynomial $2\Lambda^2\{H,F\}$
$$2\Lambda^2\{H,F\}=l_{n+1}p_1^{n+1}+l_np_1^np_2+\dots+\l_0p_2^{n+1},$$
$$
 l_{n+1}=2a_{0_x}\Lambda+na_{0}\Lambda_x+a_1\Lambda_y,
$$
$$
 l_{n}=2a_{1_x}\Lambda+(n-1)a_{1}\Lambda_x+2a_{0_y}\Lambda+2a_2\Lambda_y,
$$
$$
 l_{n-1}=2a_{2_x}\Lambda+na_{0}\Lambda_x+(n-2)a_2\Lambda_x+2a_{1_y}\Lambda+a_1\Lambda_y+3a_3\Lambda_y,
$$
$$
 l_{n-2}=2a_{3_x}\Lambda+(n-1)a_{1}\Lambda_x+(n-3)a_3\Lambda_x+2a_{2_y}\Lambda+2a_2\Lambda_y+4a_4\Lambda_y,
$$
$$
 l_{n-3}=2a_{4_x}\Lambda+(n-2)a_{2}\Lambda_x+(n-4)a_4\Lambda_x+2a_{3_y}\Lambda+3a_3\Lambda_y+5a_5\Lambda_y,
$$
$$
\dots\qquad\dots\qquad\dots\qquad\dots\qquad\dots
$$
$$
 l_3=2a_{n-2_x}\Lambda+4a_{n-4}\Lambda_x+2a_{n-2}\Lambda_x+2a_{n-3_y}\Lambda+(n-3)a_{n-3}\Lambda_y+(n-1)a_{n-1}\Lambda_y,
$$
$$
 l_2=2a_{n-1_x}\Lambda+3a_{n-3}\Lambda_x+a_{n-1}\Lambda_x+2a_{n-2_y}\Lambda+(n-2)a_{n-2}\Lambda_y+na_{n}\Lambda_y,
$$
$$
 l_1=2a_{n_x}\Lambda+2a_{n-2}\Lambda_x+2a_{n-1_y}\Lambda+(n-1)a_{n-1}\Lambda_y,
$$
$$
 l_0=a_{n-1}\Lambda_x+2a_{n_y}\Lambda+na_{n}\Lambda_y.
$$
Let us recall that by Kolokoltsov's theorem
$$
 a_{n-1}=c_1+a_{n-3}-a_{n-5}+\dots,\quad a_n=c_2+a_{n-2}-a_{n-4}+\dots
$$
We have the following system of differential equations
\begin{equation}\label{eq8}
 l_{n+1}=\dots=l_0=0.
\end{equation}

\vspace{0.3cm}

\noindent{\bf Lemma 1} {\it If $n$ is even $n=2k$, then  the system
(\ref{eq8}) has the following two conservation laws
$$
 [(na_{0}-(n-2)a_{2}+(n-4)a_{4}-\dots+(-1)^{k+1}2a_{n-2})\Lambda]_x+
$$
$$
 [(-(n-2)a_{1}+(n-4)a_{3}-\dots+(-1)^{k+1}2a_{n-3}+(-1)^{k+1}(n-1)c_1)\Lambda]_y=0,
$$
$$
 [((n-2)a_{1}-(n-4)a_{3}+\dots+(-1)^{k}2a_{n-3}+(-1)^{k+1}c_1)\Lambda]_x+
$$
$$
 [(na_{0}-(n-2)a_{2}+(n-4)a_{4}-\dots+(-1)^{k+1}2a_{n-2}+(-1)^{k+1} nc_2)\Lambda]_y=0.
$$
If $n$ is odd $n=2k+1$ then  the system (\ref{eq8}) has the
following two conservation laws
$$
 [((n-1)a_{0}-(n-3)a_{2}+(n-5)a_4-\dots+(-1)^{k+1}2a_{n-3}+(-1)^{k}c_1)\Lambda]_x+
$$
$$
 [(-(n-1)a_{1}+(n-3)a_{3}-(n-5)a_{5}+\dots+(-1)^{k}2a_{n-2}+(-1)^{k} nc_2)\Lambda]_y=0,
$$
$$
 [((n-1)a_{1}-(n-3)a_{3}+\dots+(-1)^{k+1}2a_{n-2}))\Lambda]_x+
$$
$$
 [((n-1)a_{0}-(n-3)a_{2}+\dots+(-1)^{k+1}2a_{n-3}+(-1)^{k+1} (n-1)c_1)\Lambda]_y=0.
$$ }

\vspace{0.3cm}

Theorems 1 and 2 immediately follow from the Lemma 1. Indeed, if we
put $c_1=0$ in Lemma 1 we get two conservation laws of the form
(\ref{eqq3}), where for $n$ even
$$P=(na_{0}-(n-2)a_{2}+(n-4)a_{4}-\dots+(-1)^{k+1}2a_{n-2})\Lambda,$$
$$Q=(-(n-2)a_{1}+(n-4)a_{3}-\dots+(-1)^{k+1}2a_{n-3})\Lambda,$$
$$R=(-na_{0}+(n-2)a_{2}-(n-4)a_{4}+\dots+(-1)^{k}2a_{n-2}+(-1)^{k}
nc_2)\Lambda.$$ Analogously for $n$ odd one has:
$$P=((n-1)a_{1}-(n-3)a_{3}+\dots+(-1)^{k+1}2a_{n-2}))\Lambda,$$
$$Q=((n-1)a_{0}-(n-3)a_{2}+\dots+(-1)^{k+1}2a_{n-3})\Lambda,$$
$$R=(-(n-1)a_{1}+(n-3)a_{3}-(n-5)a_{5}+\dots+(-1)^{k}2a_{n-2}+(-1)^{k}
nc_2)\Lambda.$$

\vspace{0.3cm}

\noindent {\bf Proof of Lemma 1}

Let us consider the case of even degree $n=2k$. Then
$$
 a_{n-1}=c_1+a_{n-3}-a_{n-5}+\dots+(-1)^ka_1,\ a_n=c_2+a_{n-2}-a_{n-4}+\dots+(-1)^{k+1}a_0.
$$
Let us consider the following linear combination of equations (\ref{eq8})
$$
 nl_{n+1}-(n-2)l_{n-1}+(n-4)l_{n-3}-\dots+(-1)^{k+1}2l_3=
$$
$$
 2(na_{0_x}-(n-2)a_{2_x}+(n-4)a_{4_x}-\dots+(-1)^{k+1}2a_{n-2_x})\Lambda+
$$
$$
 \left(n^2a_0-(n-2)na_0-(n-2)^2a_2+(n-4)(n-2)a_2+(n-4)^2a_4-...\right.
$$
$$
 \left.+(-1)^{k+1}2\cdot4a_{n-4}+(-1)^{k+1}2^2a_{n-2}\right)\Lambda_x+
$$
$$
 2(-(n-2)a_{1_y}+(n-4)a_{3_y}-\dots+(-1)^{k+1}2a_{n-3_y})\Lambda+
$$
$$
 (n a_1-(n-2)a_1-(n-2)3a_3+(n-4)3a_3+(n-4)5a_5-\dots+(-1)^{k+1}2(n-3) a_{n-3}+
$$
$$
 (-1)^{k+1}2 (n-1)a_{n-1})\Lambda_y=$$
$$
 2(na_{0_x}-(n-2)a_{2_x}+(n-4)a_{4_x}-\dots+(-1)^{k+1}2a_{n-2_x})\Lambda+
$$
$$
 2(na_0-(n-2)a_2+(n-4)a_4-\dots+(-1)^{k+1}2 a_{n-2})\Lambda_x+
$$
$$
 2(-(n-2)a_{1_y}+(n-4)a_{3_y}-\dots+(-1)^{k+1}2a_{n-3_y})\Lambda+
$$
$$
 2(a_1-3a_3+5a_5-\dots+(-1)^{k}(n-3) a_{n-3}+(-1)^{k+1}a_{n-1}(n-1))\Lambda_y=$$
$$
 2[(na_{0}-(n-2)a_{2}+(n-4)a_{4}-\dots+(-1)^{k+1}2a_{n-2})\Lambda]_x+
$$
$$
 2(-(n-2)a_{1_y}+(n-4)a_{3_y}-\dots+(-1)^{k+1}2a_{n-3_y})\Lambda+
$$
$$
 2(a_1-3a_3+5a_5-\dots+(-1)^{k}(n-3) a_{n-3}+(-1)^{k+1}(n-1)(c_1+a_{n-3}-a_{n-5}+\dots+(-1)^ka_1))\Lambda_y=
$$
$$
 2[(na_{0}-(n-2)a_{2}+(n-4)a_{4}-\dots+(-1)^{k+1}2a_{n-2})\Lambda]_x+
$$
$$
 2(-(n-2)a_{1_y}+(n-4)a_{3_y}-\dots+(-1)^{k+1}2a_{n-3_y})\Lambda+
$$
$$
 2(-(n-2)a_1+(n-4)a_3-\dots+(-1)^{k+1}2a_{n-3})\Lambda_y+2(-1)^{k+1}(n-1)c_1\Lambda_y=$$
$$
 2[(na_{0}-(n-2)a_{2}+(n-4)a_{4}-\dots+(-1)^{k+1}2a_{n-2})\Lambda]_x+
$$
$$
 2[(-(n-2)a_{1}+(n-4)a_{3}-\dots+(-1)^{k+1}2a_{n-3})\Lambda]_y+2(-1)^{k+1}(n-1)c_1\Lambda_y=0.
$$
We have the first required conservation law. By similar calculations,
$$
 nl_{n}-(n-2)l_{n-2}+(n-4)l_{n-4}+\dots+(-1)^{k+1}2l_2=
$$
$$
 2[((n-2)a_{1}-(n-4)a_{3}+\dots+(-1)^{k}2a_{n-3}+(-1)^{k+1}c_1)\Lambda]_x+
$$
$$
 2[(na_{0}-(n-2)a_{2}+(n-4)a_{4}-\dots+(-1)^{k+1}2a_{n-2}+(-1)^{k+1} nc_2)\Lambda]_y=0.
$$
Let us consider the case of odd degree $n=2k+1$. Then
$$
 a_{n-1}=c_1+a_{n-3}-a_{n-5}+\dots+(-1)^{k+1}a_0,\ a_n=c_2+a_{n-2}-a_{n-4}+\dots+(-1)^{k+1}a_1.
$$
By direct calculation we get
$$
 (n+1)l_{n+1}-(n-1)l_{n-1}+(n-3)l_{n-3}-\dots+(-1)^{k}2l_2=
$$
$$
 2[((n-1)a_{0}-(n-3)a_{2}+(n-5)a_4-\dots+(-1)^{k+1}2a_{n-3}+(-1)^{k}c_1)\Lambda]_x+
$$
$$
 2[(-(n-1)a_{1}+(n-3)a_{3}-(n-5)a_{5}+\dots+(-1)^{k}2a_{n-2}+(-1)^{k} nc_2)\Lambda]_y=0,
$$
and
$$
 (n-1)l_{n}-(n-3)l_{n-2}+(n-5)l_{n-4}-\dots+(-1)^{k+1}2l_3=
$$
$$
 2[((n-1)a_{1}-(n-3)a_{3}+\dots+(-1)^{k+1}2a_{n-2}))\Lambda]_x+
$$
$$
 2[((n-1)a_{0}-(n-3)a_{2}+\dots+(-1)^{k+1}2a_{n-3}+(-1)^{k+1} (n-1)c_1)\Lambda]_y=0.
$$
Lemma 1 and Theorems 1,2 are proved.

\section{Proof of Theorems 3 and 4}
Let us consider the case $c_1=0, n=3$ in Lemma 1. Then we have
conservation laws (\ref{eqq3}) where
\begin{equation}\label{eq9}
 P=-a_1\Lambda,\quad Q=-a_0\Lambda,\quad R=\left(\frac{3}{2} c_2+a_1\right)\Lambda.
\end{equation}
From (\ref{eqq3}) it follows that
\begin{equation}\label{eq10}
 P=h_{yy},\quad Q=-h_{xy},\quad R=h_{xx},
\end{equation}
where $h(x,y)$ is some function. From (\ref{eq9}) and (\ref{eq10}) we obtain
\begin{equation}\label{eq11}
 a_0=\frac{3c_2h_{xy}}{2\Delta h},\quad a_1=-\frac{3c_2h_{yy}}{2\Delta h},\quad \Lambda=\frac{2\Delta h}{3c_2}.
\end{equation}
For $n=3$ the equation $l_4=0$ has the form
$$
 a_1\Lambda_y+2a_{0_x}\Lambda+3a_0\Lambda_x=0.
$$
Let us substitute (\ref{eq11}) in the last equation. We get
$$
 2h_{xx}h_{xxy}+h_{yy}(h_{xxy}-h_{yyy})+h_{xy}(h_{xxx}+h_{xyy})=0.
$$
Since $a_0,a_1,\Lambda$ are periodic function we have
$$
 h=\lambda+a_{11}x^2+2a_{12}xy+a_{22}y^2,$$
 where $\lambda$ is a function periodic with respect to $\Gamma$.
Then $\lambda$ satisfies the equation (\ref{eq4}). Theorem 3 is
proved.

Let us consider the case $c_1=0, n=4$ in Lemma 1. Then we have
conservation laws (\ref{eqq3}) where
\begin{equation}\label{B1}
 P=(2a_0-a_2)\Lambda,\quad Q=-a_1\Lambda,\quad R=(2c_2-2a_0+a_2)\Lambda.
\end{equation}
We have
\begin{equation}\label{B2}
 P=f_{yy},\quad Q=-f_{xy},\quad R=f_{xx},
\end{equation}
where $f(x,y)$ is some function.
From (\ref{B1}), (\ref{B2}) we get
\begin{equation}\label{B3}
 a_1=2c_2\frac{f_{xy}}{\Delta f},\quad a_2=-2c_2\frac{f_{yy}}{\Delta f}+2a_0,\quad \Lambda=\frac{\Delta f}{2c_2}.
\end{equation}
Using (\ref{B3}) from $l_5=0$ and $l_4=0$ we get
\begin{equation}\label{B4}
 a_{0_x}=-\frac{1}{(\Delta f)^2}(c_2f_{yyy}f_{xy}+c_2f_{xy}f_{xxy}+2a_0\Delta f(f_{xxx}+f_{xyy})),
\end{equation}
$$
 a_{0_y}=\frac{1}{(\Delta f)^2}(2f_{yy}((c_2-a_0)f_{yyy}-a_0f_{xxy})
$$
\begin{equation}\label{B5}
 -2f_{xx}(c_2f_{xxy}+a_0(f_{yyy}+f_{xxy}))-c_2f_{xy}(f_{xyy}+f_{xxx})).
\end{equation}
We differentiate (\ref{B4}) with respect to $y$, (\ref{B5}) --- with respect to $x$ and take a difference between the results, and after that
we substitute into the result instead of $a_{0_x}$ and $a_{0_y}$ expressions (\ref{B4}) and (\ref{B5}). It gives  an equation on $f$
$$
 f_{xy}(f_{yyyy}-f_{xxxx})+3(f_{yyy}f_{xyy}-f_{xxx}f_{xxy})+2(f_{yy}f_{xyyy}-f_{xx}f_{xxxy})=0.
$$
Since $P,Q$ and $R$ are periodic functions, $f$ can be written in the form
$$
 f=\lambda+a_{11}x^2+2a_{12}xy+a_{22}y^2,
$$
where $\lambda$ is a function periodic with respect to $\Gamma$.
This yields (\ref{e4}). Theorem 4 is proved.

\section{Proof of Theorems 5 and 6}
Let us prove first Theorem 5. We start with the simple case of
$n=3$. Given a periodic solution $\lambda$ of equation (\ref{eq4})
satisfying $\Delta\lambda+2(a_{11}+a_{22})>0,$ then it follows from
the explicit formulas (\ref{eq11}) that the coefficients of the
integral $a_0, a_1, a_2 , a_3$ as well as the factor $\Lambda$ are
periodic functions.

In order to treat the case $n=4$ we proceed with a topological argument as follows.
Let $\lambda$ be a periodic solution of the equation (\ref{e4})
satisfying $\Delta\lambda+2(a_{11}+a_{22})>0.$ Coefficients of the integral and the conformal factor
$\Lambda$ are determined by the function $$
 f=\lambda+a_{11}x^2+2a_{12}xy+a_{22}y^2
$$
using the formulas (\ref{B3}), (\ref{B4}), (\ref{B5}). Notice that by Kolokoltsov identities and
(\ref{B3}) the functions $\Lambda, a_1, a_3$ are periodic with respect to the lattice.
However, coefficients  $a_0$ and also $a_2, a_4$ are not necessarily periodic. Coefficient $a_0$
is determined by the equations (\ref{B4}), (\ref{B5}) and it is convenient to rewrite them in the form

\begin{equation}\label{B4'}
(a_0(\Delta f)^2)_x=-c_2f_{xy}(\Delta f)_y:=V
\end{equation}
\begin{equation}\label{B5'}
(a_0(\Delta f)^2)_y=c_2(f_{yy}^2-f_{xx}^2)_y-c_2f_{xy}(\Delta f)_x:=W.
\end{equation}

The functions $\Delta f, V, W $ are periodic and we need to show that the periods of the $1$-form
$Vdx+Wdy$ on the torus are zeroes. We have from (\ref{B4'}), (\ref{B5'})
$$
 a_0=\frac{L(x,y)}{(\Delta f)^2}+\frac{P_0(x,y)}{(\Delta f)^2},
$$
where $L$ is a linear function and $P_0$ and $\Delta f$ are
periodic functions with respect to the lattice. Due to (\ref{B3}) we have for other coefficients an
analogous form:
$$a_2=\frac{2L(x,y)}{(\Delta f)^2}+\frac{P_2(x,y)}{(\Delta f)^2},$$
and since $a_4=c_2+a_2-a_0$ we get also
$$a_4=\frac{L(x,y)}{(\Delta f)^2}+\frac{P_4(x,y)}{(\Delta f)^2},$$
where $P_2, P_4$ are periodic functions.
In addition, odd coefficients $a_1,a_3$ are periodic functions. Take
now two periodic geodesics $\gamma_1,\gamma_2$ on the covering plane
of the configuration torus representing two independent homotopy
classes $e_1, e_2$ of the lattice. Denote by $z$ the intersection
point and by
$$
 z_1=z+e_1,\ z_2=z+e_2
$$
the translations of $z$. Since
$$
 F=a_0(x,y)p_1^4+a_1(x,y)p_1^{3}p_2+a_2(x,y)p_1^{2}p_2^2+a_3(x,y)p_1p_2^3+
 a_4(x,y)p_2^4,
$$
is the first integral of the geodesic flow we have that the
increment of $F$ along the two geodesics $\gamma_1,\gamma_2$ must
vanish. But on the other hand we compute:
$$\Delta F|_{\gamma_1}=\frac{L(e_1)}{(\Delta f)^2}(p^{'}_1)^4+\frac{2L(e_1)}{(\Delta
f)^2}(p^{'}_1)^2(p^{'}_2)^2+\frac{L(e_1)}{(\Delta f)^2}(p^{'}_2)^4,$$
$$\Delta F|_{\gamma_2}=\frac{L(e_2)}{(\Delta f)^2}(p^{''}_1)^4+\frac{2L(e_2)}{(\Delta
f)^2}(p^{''}_1)^2(p^{''}_2)^2+\frac{L(e_2)}{(\Delta f)^2}(p^{''}_2)^4,$$ where we used
the form of the coefficients $a_i$. Here we used $(p^{'}_1,p^{'}_2)$ (respectively $(p^{''}_1,p^{''}_2)$) for the momenta variables corresponding to the tangent vector $\dot{\gamma_1}(z)$ (respectively $\dot{\gamma_2}(z)$) at the intersection point $z$. But the last two identities
reduce to
$$\Delta F|_{\gamma_1}=\frac{L(e_1)}{(\Delta f)^2}((p^{'}_1)^2+(p^{'}_2)^2)^2=0,$$

$$\Delta F|_{\gamma_2}=\frac{L(e_2)}{(\Delta f)^2}((p^{''}_1)^2+(p^{''}_2)^2)^2=0.$$
Thus $$L(e_1)=L(e_2)=0,$$ which means that the linear function $L$
vanishes. This completes the proof of Theorem 5 for $n=4$.

Let us finish this section establishing variational form of the
equation (\ref{e4}). This becomes clear if one rewrites (\ref{e4})
in the following way:
$$
4a_{22}\lambda_{xyyy}-4a_{11}\lambda_{xxxy}+
2a_{12}(\lambda_{yyyy}-\lambda_{xxxx})+$$
\begin{equation}\label{e4'}
+(\lambda_{xy}\lambda_{yy})_{yy}-(\lambda_{xy}\lambda_{xx})_{xx}+
\frac{1}{2}(\lambda_{yy}^2-\lambda_{xx}^2)_{xy}=0.
\end{equation}
It is easy to verify that last equation is indeed Euler-Lagrange
equation of the functional of Theorem 6. This completes the proof.

\section{Proof of Theorem 7}
As in the proof of Theorem 6 let us rewrite equation (\ref{e5}) in
the form
$$\lambda_{xxxy}-\lambda_{yyyx}=$$
\begin{equation}\label{e4''}
=\varepsilon
(\lambda_{xy}\lambda_{yy})_{yy}-(\lambda_{xy}\lambda_{xx})_{xx}+
\frac{1}{2}(\lambda_{yy}^2-\lambda_{xx}^2)_{xy}.
\end{equation}
Then the recursion step (\ref{e6}) looks as follows:
\begin{equation}\label{e4'''}
(\lambda_k)_{xxxy}-(\lambda_k)_{yyyx}=\sum_{p,q\geq0,p+q=k-1}^{k-1}A_{pq}+B_{pq}+C_{pq},
\end{equation}
where$$ A_{pq}= ((\lambda_p)_{xy}(\lambda_q)_{yy})_{yy},$$

$$B_{pq}=-((\lambda_p)_{xy}(\lambda_q)_{xx})_{xx},$$

$$C_{pq}=\frac{1}{2}((\lambda_p)_{yy}(\lambda_q)_{yy}-(\lambda_p)_{xx}(\lambda_q)_{xx})_{xy}.$$

We prove by induction the following claim. All  $\lambda_k$ have no
Fourier monomials of the form
\begin{equation}\label{mon}
 e^{inx},\quad e^{iny},\quad e^{in(x+y)},\quad e^{in(x-y)},
\end{equation}
and is a function which is symmetric with respect to axes and
diagonals.

Assume inductively that for all $k=1,...,K-1$ the claim holds. In
order to construct $\lambda_K$ one needs that the monomials
(\ref{mon}) do not show up in the right hand side of (\ref{e4'''}).
Start with $e^{inx}$. Such a monomial can appear on the right hand
side only from $B_{pq}$. But the functions $\lambda_p, \lambda_q$
are even with respect to $y$, therefore $B_{pq}$ is odd with respect
to $y$ and so the Fourier coefficient of $e^{inx}$ must vanish.
Analogously, $e^{iny}$ can appear only from $A_{pq}.$ But this is
again an odd function on $x$ and thus the Fourier coefficient of
$e^{iny}$ must vanish. In order to conclude about the monomials
$e^{in(x+y)}$,$e^{in(x-y)}$ we notice that the equation (\ref{e4''})
and therefore also (\ref{e4'''}) is invariant on the rotation of the
plane by $45^{\circ}$ and so the previous argument can be applied.
Thus monomials (\ref{mon}) do not appear and $\lambda_k$ can be
found. One can easily see it is also symmetric with respect to the
axes and diagonals. This proves the claim.

\vspace{0.4cm}

\noindent{\bf Acknowledgments}

The first author (MB) was partially supported by ISF grant 128/10,
second author (AM) was partially supported by RFBR grant 12-01-00124-a, by
grant from Dmitri Zimin's "Dynasty" foundation and by grant of the
Russian Federation for the State Support of Researches (Agreement No
14.B25.31.0029). This paper was initiated on summer 2013 when both
authors participated in Integrable Systems Programm in CIRM
(Marseille) and CRM (Barcelona). It is a pleasure to thank the
organizers and the Institutions for their support.


\begin{thebibliography}{}

\bibitem {AA} S.V. Agapov, D.N. Alexandrov.
Fourth-Degree Polynomial Integrals of a Natural Mechanical System on
a Two-Dimensional Torus // Math. Notes. 2013. V. 93. N. 5. P.
780–-783.

\bibitem {b}
M. Bialy.
Richness or semi-hamiltonicity of quasi-linear systems that are not in evolution
form // Quarterly of Applied Math. 2013. V. 71. P. 787--796.

\bibitem{BM}
M. Bialy, A. Mironov. Rich quasi-linear system for integrable
geodesic flows on 2-torus // Discrete and Continuous Dynamical
Systems - Series A. 2011. V. 29. N. 1. P. 81--90.

\bibitem{BM1}
M. Bialy, A. Mironov. Qubic and Quartic integrals for geodesic flow
on 2-torus via system of Hydrodynamic type // Nonlinearity. 2011. V. 24.
N. 12. P. 3541-–3554.


\bibitem{KTD}  N.V. Denisova, V.V. Kozlov, D.V. Treschev.
Remarks on polynomial integrals of higher degrees for reversible
systems with toral configuration space // Izv. Math.  2012. V.
76. N. 5. P. 907-–921.


\bibitem{DN}

B.A. Dubrovin, S.P. Novikov.  Hydrodynamics of weakly deformed soliton lattices. Differential geometry and Hamiltonian theory //
Russian Math. Surveys. 1989. V. 44. N. 6. P. 35–124.

\bibitem{MF}

O.I. Mokhov, E.V. Ferapontov. Non-local Hamiltonian operators of hydrodynamic type related to metrics of constant curvature //
Russian Math. Surveys. 1990. V. 45. N. 3. P. 218–-219.


\bibitem{KD}

V.V. Kozlov, N.V. Denisova. Polynomial integrals of geodesic flows on a two-dimensional torus // Sb. Math. 1995.
V. 83. N. 4. P. 69-–81.

\bibitem{KT}

V.V. Kozlov, D.V. Treschev. On the integrability of Hamiltonian systems with toral position space // Math.
USSR-Sb. 1989. V. 63. N. 1. P. 121-–139.

\bibitem{K}

V.N. Kolokoltsov. Geodesic flows on two-dimensional manifolds with
an additional first integral that is polynomial in the velocities //
Math. USSR-Izv. 1983. V. 21. N. 2. P. 291–-306

\bibitem{PT}

M.V. Pavlov, S.P. Tsarev. Tri-Hamiltonian structures of Egorov systems of hydrodynamic type // Funct. Anal. Appl.
2003. V. 37. N. 1. P. 32-–45.

\bibitem{S}

D. Serre, "Systems of Conservation Laws," Vol. 2, Geometric structures, oscillations, and
initial-boundary value problems. Translated from the 1996 French original by I. N. Sneddon,
Cambridge University Press, Cambridge, 2000.


\bibitem{Ts}

S.P. Tsarev. The geometry of Hamiltonian systems of hydrodynamic type. The generalized hodograph method // Math.
USSR-Izv. 1991. V. 37. N. 2. P. 397–-419.




\end{thebibliography}
\end{document}